\documentclass[12pt, a4paper]{article}

\skip\footins=8mm plus 1mm

\usepackage{latexsym,amsmath,amsfonts,amssymb,amsthm,enumitem}

\long\def\delete#1{}

\usepackage{longtable, booktabs}
\usepackage{supertabular}
\usepackage{rotating}
\usepackage{multicol}
\usepackage{multirow}
\usepackage{makeidx} 
\usepackage{epsfig}

\usepackage{hyperref}
\usepackage[labelsep=period]{caption}

\newtheorem{theorem}{Theorem}[section]

\newtheorem{lemma}[theorem]{Lemma}

\newtheorem{definition}[theorem]{Definition}

\newtheorem{problem}[theorem]{Problem}

\def\proof{\par\noindent{\textbf{Proof.}~}}

\def\la{\langle}
\def\ra{\rangle}

\def\val{{\rm val}}

\def\ord{{\rm ord}}
\def\G1{G^\mathcal{C}}
\def\magma{{\sc Magma} }
\def\row{{\rm row}}

\newcommand{\bmat}[1]{\begin{bmatrix}#1\end{bmatrix}}

\newcommand{\beq}{\begin{equation}}
\newcommand{\eeq}{\end{equation}}
\newcommand{\bea}{\begin{eqnarray}}
\newcommand{\eea}{\end{eqnarray}}
\newcommand{\bean}{\begin{eqnarray*}}
\newcommand{\eean}{\end{eqnarray*}}

\def\BB{\mathcal{B}}  \def\DD{\mathcal{D}}       \def\LL{\mathcal{L}}    \def\PP{\mathcal{P}}      

  \def\lcm{{\rm lcm}}  \def\row{{\rm row}} 

  \def\Ga{\Gamma}      \def\Om{\Omega}    

\def\a{\alpha}  \def\b{\beta}  \def\d{\delta}  \def\g{\gamma}  \def\om{\omega}  \def\l{\lambda}  \def\s{\sigma}
\def\t{\tau}    \def\vp{\varphi}  

  \def\Sym{{\rm Sym}}     
\def\Aut{{\rm Aut}}    \def\soc{{\rm soc}}    
  \def\PSL{{\rm PSL}}  \def\PGL{{\rm PGL}}  \def\PSU{{\rm PSU}}  \def\PGU{{\rm PGU}}  \def\SL{{\rm SL}}  \def\GL{{\rm GL}}  \def\GU{{\rm GU}}    \def\AGL{{\rm AGL}}  \def\ASL{{\rm ASL}}  \def\PGL{{\rm PGL}}  \def\AG{{\rm AG}}  \def\PG{{\rm PG}}      \def\GammaL{{\rm\Gamma L}}
\def\AGammaL{{\rm A\Gamma L}}  \def\PGammaL{{\rm P\Gamma L}}    
  \def\Sp{{\rm Sp}}      \def\C{{\rm C}}
     \def\G{{\rm G}}  \def\R{{\rm R}}  \def\Sz{{\rm Sz}}          \def\Co{{\rm Co}}             \def\HS{{\rm HS}}          \def\GSp{{\rm GSp}}

\def\bzero{\mathbf{0}}   \def\bb{\mathbf{b}} \def\bc{\mathbf{c}}  \def\bfe{\mathbf{e}}
    \def\bx{\mathbf{x}} \def\by{\mathbf{y}}
\def\bz{\mathbf{z}}    \def\bv{\mathbf{v}} 

\usepackage{xcolor}
\usepackage[normalem]{ulem}

\begin{document}
\openup 0.5\jot

\title{A family of symmetric graphs in relation to 2-point-transitive linear spaces}

\author{\renewcommand{\thefootnote}{\arabic{footnote}}Teng Fang\footnotemark[1] , Sanming Zhou\footnotemark[2] , Shenglin Zhou\footnotemark[3]}

\footnotetext[1]{School of Mathematics, Southeast University, Nanjing 211189, P. R. China}

\footnotetext[2]{{School} of Mathematics and Statistics, The University of Melbourne, Parkville, VIC 3010, Australia}

\footnotetext[3]{School of Mathematics, South China University of Technology, Guangzhou, Guangdong, 510640, P. R. China}

\renewcommand{\thefootnote}{}
\footnotetext{{\em E--mail addresses}: \texttt{tfang@seu.edu.cn} (Teng Fang), \texttt{sanming@unimelb.edu.au} (Sanming Zhou), \texttt{slzhou@scut.edu.cn} (Shenglin Zhou)}

\date{}

\maketitle

\begin{abstract}
A graph $\Ga$ is $G$-symmetric if it admits $G$ as a group of automorphisms acting transitively on the set of arcs of $\Ga$, where an arc is an ordered pair of adjacent vertices. Let $\Ga$ be a $G$-symmetric graph such that its vertex set admits a nontrivial $G$-invariant partition $\BB$, and let $\DD(\Ga, \BB)$ be the incidence structure with point set $\BB$ and blocks $\{B\} \cup \Ga_{\BB}(\a)$, for $B \in \BB$ and $\a \in B$, where $\Ga_{\BB}(\a)$ is the set of blocks of $\BB$ containing at least one neighbour of $\a$ in $\Ga$. In this paper we classify all $G$-symmetric graphs $\Ga$ such that $\Ga_{\BB}(\a) \ne \Ga_{\BB}(\b)$ for distinct $\a, \b \in B$, the quotient graph of $\Ga$ with respect to $\BB$ is a complete graph, and $\DD(\Ga, \BB)$ is isomorphic to the complement of a $(G, 2)$-point-transitive linear space.


\end{abstract}

\section{Introduction}
\label{sec:introd}

Since Tutte's seminal work \cite{Tutte}, the study of symmetric graphs and highly arc-transitive graphs has long been an important area of research in algebraic graph theory (see \cite{Praeger97, Praeger00}). It is desirable but often challenging to classify symmetric graphs with certain combinatorial and/or group-theoretical properties. In this paper we prove a result of this kind by classifying all symmetric graphs such that a certain associated block design is a $2$-point-transitive linear space. We are able to achieve this classification owing to the classification of $2$-point-transitive linear spaces \cite{Kantor85} and the flag graph construction \cite{Zhou-flag}, the former being a consequence of the classification of finite simple groups.

All graphs considered are finite and undirected with no loops or parallel edges, and all groups considered are finite. The \emph{order} of a graph is its number of vertices, and the \emph{valency} of a regular graph is the valency (degree) of every vertex in the graph. Let $\Ga$ be a graph with vertex set $V(\Ga)$. Let $G$ be a group acting on $V(\Ga)$ as a group of automorphisms of $\Ga$ (that is, $G$ preserves the adjacency and non-adjacnecy relations of $\Ga$). If $G$ is transitive on $V(\Ga)$ and transitive on the set of arcs of $\Ga$, then $\Ga$ is said to be {\em $G$-symmetric} (or {\em $G$-arc transitive}), where an {\em arc} is an ordered pair of adjacent vertices. If in addition $V(\Ga)$ admits a nontrivial {\em $G$-invariant partition} $\BB = \{B, C, \ldots\}$, that is, $1 < |B| < |V(\Ga)|$ and $B^g := \{\a^g\;|\; \a \in B\} \in \BB$  for any $B \in \BB$ and $g \in G$, then $\Ga$ is called an \emph{imprimitive} $G$-symmetric graph and $(\Ga, G, \BB)$ is said to be a {\em symmetric triple}. In this case the {\em quotient graph} of $\Ga$ relative to $\BB$, denoted by $\Ga_{\BB}$, is defined to be the graph with vertex set $\BB$ in which $B, C \in\BB$ are adjacent if and only if there exists at least one edge of $\Ga$ between $B$ and $C$. We assume without mentioning explicitly that $\Ga_{\BB}$ has at least one edge, so that each block of $\BB$ is an independent set of $\Ga$. For $\a \in V(\Ga)$, denote by $\Ga(\a)$ the neighbourhood of $\a$ in $\Ga$, and define $\Ga_{\BB}(\a)$ to be the set of blocks of $\BB$ containing at least one neighbour of $\a$ in $\Ga$. Define $\Ga(B) := \cup_{\a \in B}\Ga(\a)$  and denote by $\Ga_{\BB}(B)$ the neighbourhood of $B$ in $\Ga_{\BB}$. Note that $\Ga_{\BB}(\a) \subseteq \Ga_{\BB}(B)$ for $\a \in B$. For adjacent blocks $B, C$ of $\BB$, define $\Ga[B, C]$ to be the induced bipartite subgraph of $\Ga$ with bipartition $\{\Ga(C) \cap B, \Ga(B) \cap C\}$. Since $\Ga$ is $G$-symmetric and $\BB$ is $G$-invariant, up to isomorphism, $\Ga[B, C]$ is independent of the choice of adjacent blocks $B, C$ of $\BB$. Define $\DD(B) := (B, \Ga_{\BB}(B))$ to be the incidence structure with ``point set'' $B$ and ``block set'' $\Ga_{\BB}(B)$ such that $\a \in B$ and $C \in \Ga_{\BB}(B)$ are incident if and only if $C \in \Ga_{\BB}(\a)$ (see \cite{Gardiner-Praeger95}). Set
$$
v := |B|,\;\; r := |\Ga_{\BB}(\a)|,\;\; b := |\Ga_{\BB}(B)|,\;\; k := |\Ga(C) \cap B|.
$$
These parameters depend on $(\Ga, \BB)$ but are independent of the choice of $\a \in V(\Ga)$ and adjacent $B, C \in \BB$. It can be verified \cite{Gardiner-Praeger95} that $\DD(B)$ is a $1$-$(v, k, r)$ design with $b$ blocks, and is independent of $B$ up to isomorphism. A geometrical approach introduced by Gardiner and Praeger \cite{Gardiner-Praeger95} suggests that one may study imprimitive symmetric graphs $\Ga$ through the analysis of $\Ga_{\BB}$, $\Ga[B, C]$ and $\DD(B)$. See \cite{FFXZ, Gardiner-Praeger18, Li-Praeger-Zhou98, Zhou98, Zhou-EJC, Zhou-flag} for several studies of imprimitive symmetric graphs using this approach and \cite{CZ14, FFXZ2, SZ13, Zhou-EJC} for a complete classification of symmetric triples $(\Ga, G, \BB)$ such that $k=v-1 \ge 2$ and $\Ga_{\BB}$ is a complete graph.

As in \cite{Zhou-flag}, let $\DD^{*}(B) := (\Ga_{\BB}(B), B)$ be the dual $1$-design of $\DD(B)$. Here we may identify each ``block'' $\a \in B$ of $\DD^{*}(B)$ with the subset $\Ga_{\BB}(\a)$ of the ``point set'' $\Ga_{\BB}(B)$ of $\DD^{*}(B)$, and we call two such blocks $\Ga_{\BB}(\b)$, $\Ga_{\BB}(\g)$ (for distinct $\b, \g \in B$) {\em repeated} if $\Ga_{\BB}(\b) = \Ga_{\BB}(\g)$. Define $\DD(\Ga, \BB)$ to be the 1-design with point set $\BB$ and block set
$$
(\{B\} \cup \Ga_{\BB}(\a))^G = \{\{B(\d)\} \cup \Ga_{\BB}(\d): \d \in V(\Ga)\},
$$
with repeated blocks identified, where $B(\d)$ is the unique block of $\BB$ containing $\d$. Then $\DD(\Ga, \BB)$ admits $G$ as a point- and block-transitive group of automorphisms (see \cite[Lemma 3.1]{Zhou-flag}), and $\DD^{*}(B)$ admits $G_B$ (the setwise stabilizer of $B$ in $G$) as a point- and block-transitive group of automorphisms. It would be interesting to investigate how $\DD(\Ga, \BB)$ and $\Ga$ interact with each other and to what extent the structure of the former determines that of the latter. In this regard the case where $\Ga_{\BB}$ is a complete graph (hence $\Ga_{\BB}\cong K_{b+1}$) is most interesting because in this case $\DD(\Ga, \BB)$ is a $(G,2)$-point-transitive $2$-$(b+1, r+1, \l)$ design for some $\l \geq 1$ (see \cite{CZ14, FFXZ, FFXZ2, Gardiner-Praeger18, SZ13} for related works). In particular, the following general problem may be worth studying (see \cite[Section 3]{Zhou-flag} for related discussions).

\begin{problem}
Classify all symmetric triples $(\Ga, G, \BB)$ such that $\Ga_{\BB}$ is a complete graph and $\DD(\Ga, \BB)$ belongs to a given family of $(G,2)$-point-transitive $2$-designs.
\end{problem}

In this paper we address this problem in the case where (i) $\Ga_{\BB}\cong K_{b+1}$, (ii) $\DD^*(B)$ contains no repeated blocks, and (iii) the complement of $\DD(\Ga, \BB)$ is a linear space, and we give a classification of all such $G$-symmetric graphs $\Ga$ in this case. Here the third condition is equivalent to saying that for any distinct blocks $C_1, C_2 \in \BB$ and any vertices $\a, \b \in V(\Ga)$ such that $\a, \b \not \in C_i$ and $\Ga(\a) \cap C_i = \Ga(\b) \cap C_i = \emptyset$ for $i=1,2$, we have $\{B(\a)\} \cup \Ga_{\BB}(\a) = \{B(\b)\} \cup \Ga_{\BB}(\b)$. In general, a {\em linear space} \cite{Beth-Jung-Lenz} is an incidence structure of points and lines such that every line contains at least two points and every pair of points is contained in exactly one line. A linear space is \emph{proper} if every line contains at least three points, and {\em improper} otherwise. Thus an improper 2-point-transitive linear space with $u$ points can be identified with the complete graph $K_{u}$ of order $u$. The main result in this paper is as follows.

\delete
{
Suppose $\Ga_{\BB}\cong K_{b+1}$ is complete. Then $\DD(\Ga, \BB)$ is a $(G,2)$-point-transitive $2$-$(b+1, r+1, \l)$ design, with $\frac{\l b(b+1)}{r(r+1)}$ blocks and with $\frac{\l b}{r}$ blocks per point. The complement of $\DD(\Ga, \BB)$ is a $2$-$\left(b+1, b-r, \frac{\l(b-r)(b-r-1)}{r(r+1)}\right)$ design, with $\frac{\l b(b+1)}{r(r+1)}$ blocks and with $\frac{\l b(b-r)}{r(r+1)}$ blocks per point. The complement of $\DD(\Ga, \BB)$ is a linear space if and only if $\l = \frac{r(r+1)}{(b-r)(b-r-1)}$, or equivalently $\DD(\Ga, \BB)$ has exactly $\frac{b(b+1)}{(b-r)(b-r-1)}$ blocks. Since the number of blocks of $\DD(\Ga, \BB)$ with multiplicity counted is equal to $v(b+1) = \frac{bk(b+1)}{r}$, each block of $\DD(\Ga, \BB)$ is repeated $v(b+1)/\{\frac{\l b(b+1)}{r(r+1)}\} = \frac{(r+1)k}{\l}$ times. In particular, $\l$ is a divisor of $(r+1)k$. So the complement of $\DD(\Ga, \BB)$ is a linear space if and only if for any two distinct $C_1, C_2 \in \BB$ there are exactly $\frac{(r+1)k}{\l}$ vertices $\a \in V(\Ga)$ such that $\a \not \in C_i$ and $\Ga(\a) \cap C_i = \emptyset$ for $i=1,2$.
}

\begin{theorem}
\label{thm:main theorem}
Let $\Ga$ be a $G$-symmetric graph such that $V(\Ga)$ admits a nontrivial $G$-invariant partition $\BB$, where $G \le \Aut(\Ga)$. Suppose that $\Ga_{\BB}$ is a complete graph, $\DD^{*}(B)$ contains no repeated blocks, and $\DD(\Ga, \BB)$ is isomorphic to the complement of a $(G,2)$-point-transitive linear space $\DD$.
\begin{itemize}
\item[\rm (a)] If $\DD$ is proper, then $\DD= \PG(2,q)$, $\PSL(3,q)\leq G\leq\PGammaL(3,q)$, $\Ga$ has order $q^2(q^2+q+1)$ and valency $\val(\Ga) \in \{q^2-1, (q^2-1)(q-1), (q^2-1)(q-1)m/3, (q^2-1)q\cdot\lcm(\rho, s)/s\}$, and for any distinct blocks $B, C \in \BB$ there are exactly $q(q-1)$ vertices in $B$ that have neighbours in $C$, where $m$, $s$ and $\rho$ are integers determined by $\Ga$ and $G$.
\item[\rm (b)] If $\DD$ is improper, then $\DD$, $G$ and $\Ga$ are given in Table \ref{tab0a}, where $u = |\BB|$.
\end{itemize}
\end{theorem}

\begin{table}[ht]
\begin{center}
\scalebox{0.8}[0.9]{
  \begin{tabular}{l|l|l|l|l}
\hline
     & $G$  & $\DD$ & $\Ga$  & Details \\ \hline
(b1) & $\soc(G)=A_u$, & $2$-$(u, 2, 1)$ & $\ord={u(u-1)(u-2)}/{2}$; & \textsection\ref{sec:almost simple 2-transitive groups}, case (i) \\
     & $u\geq5$ &   & $\val=u-3$ if $G=A_5$; &  \\
     &  &   & $\val=u-3$ or $2(u-3)$ if $G=S_5$ or $A_6$; & \\
     &  &   & $\val=u-3$ or $4(u-3)$ if $G=S_6$; & \\
     & &   &  $\val=u-3$, $(u-3)(2u-8)$ or $\frac{(u-3)(u-4)(u-5)}{2}$ & \\
     & & & if $u>6$ & \\\hline

(b2) & $\PGammaL(2,8)$ & $2$-$(9, 2, 1)$ & $\ord=252$ and $\val=6$ & \textsection\ref{sec:almost simple 2-transitive groups}, case (ii) \\\hline

(b3) & $A_7$ & $2$-$(15, 2, 1)$  & $\ord=105$, $\val=12$ or $24$ & \textsection\ref{sec:almost simple 2-transitive groups}, case (xi) \\\hline

(b4) & $M_u$, $u=11, 12$, & $2$-$(u, 2, 1)$ & $\ord=495$, $\val=8$ or $16$ if $u=11$; & \textsection\ref{sec:almost simple 2-transitive groups}, case (xii) \\
    & $23$ or $24$  &   & $\ord=660$, $\val=9$, $36$, $72$ or $144$ if $u=12$; &   \\
    &   &   & $\ord=5313$, $\val=20$, $60$, $120$ or $480$ if $u=23$; & \\
    &   &   & $\ord=6072$, $\val=21$, $630$, $840$ or $3360$ if $u=24$ & \\\hline

(b5) & $\AGL(1,4)$ or & $2$-$(4, 2, 1)$  & $\ord=12$ and $\val=1$ & \textsection\ref{sec:affine simple 2-transitive groups}, case (i) \\
      & $\AGammaL(1,4)$ &    &   &   \\\cline{2-4}
    & $\AGL(1,5)$ & $2$-$(5, 2, 1)$  & $\ord=10$ and $\val=2$ &   \\ \hline

(b6) & $\mathbb{F}_2^2\rtimes\la\bmat{1 & 1\\1 & 0}\ra$ or & $2$-$(4, 2, 1)$  & $\ord=12$ and $\val=1$ & \textsection\ref{sec:affine simple 2-transitive groups}, case (ii) \\
      & $\AGL(2,2)$ &   &  &  \\\hline

(b7) & $\ASL(2,3)$ & $2$-$(9, 2, 1)$  & $\ord=36$ and $\val=6$ & \textsection\ref{sec:affine simple 2-transitive groups}, case (ii) \\\cline{2-2}\cline{4-4}
   & $\AGL(2,3)$ &  & $\ord=36$, $\val=6$ or $12$ &   \\\cline{2-5}

    & $\ASL(n,3)$ or & $2$-$(3^n, 2, 1)$  & $\ord={3^n(3^n-1)}/{2}$, $\val=3^n-3$, $2(3^n-3)$ or & \textsection\ref{sec:affine simple 2-transitive groups}, case (ii) \\
    & $\AGL(n,3)$, $n\geq3$ &  & ${(3^n-3)(3^n-9)}/{2}$ & \\\hline

\end{tabular}}
\caption{Theorem \ref{thm:main theorem}: The case where $\DD$ is an improper $(G,2)$-point-transitive linear space with $|\BB|$ points. Acronym: ord = Order, val = Valency}
\label{tab0a}
  \end{center}
\end{table}

More information about the graphs in part (a)  can be found in Section \ref{sec:almost simple}, and definitions and properties of the graphs in part (b) will be given during the proofs in Sections \ref{sec:almost simple 2-transitive groups} and \ref{sec:affine simple 2-transitive groups}.

Part (a) of Theorem \ref{thm:main theorem} will be proved in Section \ref{sec:parta} using the classification \cite[Theorem 1]{Kantor85} of 2-point-transitive linear spaces. Part (b) of Theorem \ref{thm:main theorem} will be proved in Section \ref{sec:trivial} with the help of the classification of 2-transitive groups. A major tool for both parts is the flag graph construction introduced in \cite{Zhou-flag}. In a few cases \magma is used to compute the graphs involved.

\section{Preliminaries}
\label{sec:prelim}

\subsection{Notation and terminology}
\label{sec:not}

The reader is referred to \cite{Dixon-Mortimer}, \cite{Beth-Jung-Lenz} and \cite{Dembowski} for notation and terminology on permutation groups, block designs and finite geometries, respectively. Unless stated otherwise, all designs in the paper are assumed to have no repeated blocks, and we identify each block of a design with the set of points incident with it. The number of points, number of blocks, size of each block and number of blocks incident with a given point in a design $\DD$ are denoted by $v_{\DD}$, $b_{\DD}$, $k_{\DD}$ and $r_{\DD}$, respectively. In addition, if $\DD$ is a $2$-design, then the number of blocks containing a pair of points is denoted by $\l_{\DD}$.

Let $G$ be a group acting on a set $\Om$. That is, for any $\a \in \Om$ and $g \in G$ there corresponds a point in $\Om$ denoted by $\a^g$, such that $\a^{1_G} = \a$ and $(\a^g)^h = \a^{gh}$ for any $\a \in \Om$ and $g, h \in G$, where $1_G$ is the identity element of $G$. If $\a$ is a point and $P$ a subset of $\Om$, then denote by $G_{\a}$ the stabilizer of $\a$ in $G$, $G_P$ the setwise stabilizer of $P$ in $G$, and $G_{\a, P}$ the setwise stabilizer of $P$ in $G_{\a}$.

The {\em socle} of a group $G$, denoted by $\soc(G)$, is the subgroup generated by the minimal normal subgroups of $G$. It is widely known that any finite $2$-transitive group $G$ is either {\em almost simple} (that is, $\soc(G)$ is a nonabelian simple group) or {\em affine} (that is, $\soc(G)$ is elementary abelian).

The {\em complement} of an incidence structure $\DD = (\PP, \LL)$, denoted by $\overline{\DD}$, is the incidence structure with point set $\PP$ such that $L \subseteq \PP$ is a block if and only if $\PP \setminus L \in \LL$ is a block of $\DD$. In other words, $(\s, L)$ is a flag of $\overline{\DD}$ if and only if $(\s, \PP \setminus L)$ is an antiflag of $\DD$, and vice versa.

\subsection{Flag graphs}
\label{sec:flag graph}

For a symmetric triple $(\Ga, G, \BB)$, if $\DD^{*}(B)$ has no repeated blocks, then $\Ga$ can be reconstructed from $\DD(\Ga, \BB)$ by the following ``flag graph construction''.

\begin{definition}
\label{def:flag}
{\em (\cite[Definition 2.1]{Zhou-flag}) Let $\DD$ be a $G$-point-transitive and $G$-block-transitive $1$-design with block size at least $2$. A $G$-orbit $\Om$ on the set of flags of $\DD$ is said to be {\em feasible} if
\begin{itemize}
\item[(a)] for some (and hence all) point $\s$ of $\DD$, $|\Om(\s)| \ge 2$, where $\Om(\s)$ is the set of flags of $\Om$ with point-entry $\s$; and
\item[(b)] for some (and hence all) flag $(\s, L) \in \Om$, $G_{\s, L}$ is transitive on $L \setminus \{\s\}$.
\end{itemize}

Let $\Om$ be a feasible $G$-orbit on the set of flags of $\DD$, and let $\Psi$ be a $G$-orbit on
\begin{align}
\label{equ:C(D,Om)}
\C(\DD, \Om):=\{((\s,L), (\t,N))\in \Om\times\Om: \s\neq\t\text{\;and\;}\s, \t\in L\cap N\}.
\end{align}
If $\Psi$ is {\em self-paired} (that is, $((\s, L), (\t, N)) \in \Psi$ implies $((\t, N)$, $(\s, L)) \in \Psi$), then define the {\em $G$-flag graph} of $\DD$ with respect to $(\Om, \Psi)$ to be the graph with vertex set $\Om$ in which two ``vertices'' $(\s, L)$ and $(\t, N)$ are adjacent if and only if $((\s, L), (\t, N)) \in \Psi$. Denote this graph by $\Ga(\DD, \Om, \Psi)$.
}
\end{definition}


{Obviously, $\C(\DD,\Om)$ is invariant under the induced action of $G$ on $\Om\times\Om$.}

The following result was proved in \cite[Theorem 1.1]{Zhou-flag}: Given a symmetric triple $(\Ga, G, \BB)$, if $\DD^{*}(B)$ contains no repeated blocks, then $\Ga$ is isomorphic to a $G$-flag graph of $\DD(\Ga, \BB)$ with respect to some $(\Om, \Psi)$. Conversely, any $G$-flag graph $\Ga(\DD, \Om, \Psi)$ is a $G$-symmetric graph that admits $\BB(\Om) := \{\Om(\s): \s \text{\;is a point of\;} \DD\}$ as a $G$-invariant partition such that the corresponding $\DD^{*}(\Om(\s))$ contains no repeated blocks.

In the special case where $\Ga_{\BB} \cong K_{b+1}$ is a complete graph, $G$ is $2$-transitive on $\BB$ and so $\DD(\Ga, \BB)$ is a $2$-$(b+1, r+1, \l)$ design for some $\l \ge 1$ that admits $G$ as a $2$-point-transitive and block-transitive group of automorphisms. Therefore, \cite[Theorem 1.1]{Zhou-flag} mentioned above has the following corollary.

\begin{lemma}
\label{lem:two design}
{\em (\cite[Corollary 3.3]{Zhou-flag})}
Let $b \geq 2$ and $r \geq 1$ be integers, and let $G$ be a group. The following two {statements} are equivalent{:}

{\em (a)} $\Ga$ is a $G$-symmetric graph admitting a nontrivial $G$-invariant
partition $\BB$ such that $\DD^{*}(B)$ has block size $r$ and contains no
repeated blocks, and such that $\Ga_{\BB} \cong K_{b+1}${;}

{\em (b)} $\Ga \cong \Ga(\DD,\Om,\Psi)$, for a $(G,2)$-point-transitive and
$G$-block-transitive $2$-$(b+1, r+1, \l)$ design $\DD$, a feasible $G$-orbit
$\Om$ on the set of flags of $\DD$, and a self-paired $G$-orbit $\Psi$ on $\C(\DD, \Om)$.
\end{lemma}


Thus, if one wants to classify symmetric triples $(\Ga, G, \BB)$ such that $\Ga_{\BB}$ is complete, $\DD^{*}(B)$ has no repeated blocks and $\DD(\Ga, \BB)$ belongs to a given family of $(G, 2)$-point-transitive and $G$-block-transitive $2$-designs, it suffices to classify the $G$-flag graphs in Lemma \ref{lem:two design} for such designs $\DD$ in the family. We will prove Theorem \ref{thm:main theorem} in Sections \ref{sec:parta} and \ref{sec:trivial} by applying this method to the complements of $2$-point-transitive linear spaces.

\begin{lemma}
\label{lem:properties about the G-flag graph}
Let $\DD= (\PP, \LL)$ be a $(G, 2)$-point-transitive linear space. Let $\Om = (\s, \PP \setminus L)^G$ be a feasible $G$-orbit on the set of flags of $\overline{\DD}$, where $(\s, L)$ is an antiflag of $\DD$. Let $\Psi=((\s, \PP \setminus L), (\t, \PP \setminus N))^G$ be a self-paired $G$-orbit on $\C(\overline{\DD},\Om)$. Then the $G$-flag graph $\Ga:=\Ga(\overline{\DD}, \Om, \Psi)$ has order $|\PP|\cdot |L^{G_\s}|$ and valency $(|\PP|-|L|-1)\cdot|N^{G_{\s,\t, L}}|$. Moreover, there are exactly $|L^{G_{\s,\t}}|$ vertices in $\Om(\s)$ that have neighbours in $\Om(\t)$, and each of them has exactly $|N^{G_{\s,\t, L}}|$ neighbours in $\Om(\t)$.
\end{lemma}

\proof
Obviously, the number of vertices of $\Ga$ is $|\Om|=|\PP|\cdot |\Om(\s)|=|\PP|\cdot |L^{G_\s}|$. If $(\s, \PP \setminus L_1)$ has a neighbour $(\t, \PP \setminus N_1)$ in $\Om(\t)$, then $((\s, \PP \setminus L_1), (\t, \PP \setminus N_1))\in \Psi$ and thus there exists $g\in G$ such that $(\s, \PP \setminus L_1)=(\s, \PP \setminus L)^g$ and $(\t, \PP \setminus N_1)=(\t, \PP \setminus N)^g$. Hence $L_1=L^g$ and $g\in G_{\s,\t}$. Conversely, if $M$ is an image of $L$ under some element in $G_{\s,\t}$, then $(\s, \PP\setminus M)$ has a neighbour in $\Om(\t)$. Similarly, the set of neighbours of $(\s, \PP \setminus L)$ in $\Om(\t)$ is $\{(\t, \PP \setminus N^g):g\in G_{\s,\t, L}\}$ and thus the valency of $\Ga[\Om(\s), \Om(\t)]$ is $|N^{G_{\s,\t, L}}|$. Since $G_{\s,L}$ is transitive on $\PP\setminus\{L\cup\{\s\}\}$ (see Definition \ref{def:flag}), we have
$$
\{\Om(\a):\a\in\PP,\, \Ga((\s, \PP \setminus L))\cap\Om(\a)\neq\emptyset\}=\{\Om(\a):\a\in\PP\setminus\{L\cup\{\s\}\}\}.
$$
It follows that the valency of $\Ga$ is equal to $(|\PP|-|L|-1)\cdot|N^{G_{\s,\t, L}}|$.
\qed

\section{Necessary conditions}
\label{sec:some facts and necessary conditions}

\begin{lemma}
\label{lem:at most one feasible orbit}
Let $\DD$ be a $(G,2)$-point-transitive and $G$-block-transitive $2$-$(|\PP|, r+1, \l)$ design with point set $\PP$, where $r > 1$. Then there is at most one feasible $G$-orbit on the flag set of $\DD$.
\end{lemma}


\proof
Let {$\s,\t \in \PP$} be distinct points. Denote by $L_1, \ldots, L_{\lambda}$ the $\lambda$ blocks of $\DD$ containing both $\s$ and $\t$. Set $\Om_i := (\s, L_i)^G$, $i=1, \ldots, \lambda$. Since $G$ is $2$-transitive on $\PP$, we can see that $\Om_1, \ldots, \Om_{\l}$ are all possible $G$-orbits on the flag set of $\DD$ (but it may happen that $\Om_i = \Om_j$ for distinct $i$ and $j$).

Suppose that $\Om_i \ne \Om_j$ and both of them are feasible. Since $\DD$ is $G$-block-transitive, there exists a point $\xi$ of $\DD$ such that $(\xi, L_j) \in \Om_i$. The assumption $\Om_i\neq \Om_j$ implies that $\s \neq \xi$. If $G_{L_j}$ is transitive on $L_j$, then $(\eta, L_j)\in\Om_i$ for any $\eta\in L_j$ and thus $\Om_i= \Om_j$, which contradicts our assumption. Thus $G_{L_j} = G_{\xi, L_j} \le G_{\xi}$. Since $\xi \in L_j \setminus \{ \s \}$, $G_{\s, L_j} \le G_{L_j} \le G_{\xi}$ and $|L_j| = r+1 \geq 3$, $G_{\s, L_j}$ cannot be transitive on $L_j \setminus \{ \s \}$. Hence there is at most one feasible $G$-orbit on the flag set of $\DD$.
\qed

\medskip
Assume that $\DD = (\PP, \LL)$ is a $(G, 2)$-point-transitive linear space. Let $\Om = (\s, \PP \setminus L)^G = \{(\s^g, \PP \setminus L^g): g \in G\}$ be a $G$-orbit on the set of flags of $\overline{\DD}$, where $(\s, L)$ is an antiflag of $\DD$. Since $(\PP \setminus L) \setminus \{\s\} \neq \PP \setminus \{\s\}$ and $G$ is $2$-transitive on $\PP$, condition (a) in Definition \ref{def:flag} is satisfied, and hence $\Om$ is feasible if and only if $G_{\s, \PP \setminus L} = G_{\s, L}$ is transitive on $(\PP \setminus L) \setminus \{\s\}$. Moreover, if there exists a self-paired $G$-orbit on $\C(\overline{\DD},\Om)$, then we must have $|\PP|\geq4$. Hence, in the rest of this paper, we assume
$$
|\PP|\geq4.
$$
To prove Theorem \ref{thm:main theorem} it suffices to find out all feasible $\Om$ and then all self-paired $G$-orbits on $\C(\overline{\DD},\Om)$. The following result will be used to prove the nonexistence of such feasible $\Om$ for some $\DD$.


\begin{lemma}
\label{lem:two 2-transitive actions}
Let $\DD = (\PP, \LL)$ be a $(G, 2)$-point-transitive linear space, where $G\leq\Aut(\DD)$. Let $\Om = (\s, \PP \setminus L)^G$ be a feasible $G$-orbit on the set of flags of $\overline{\DD}$, where $(\s, L)$ is an antiflag of $\DD$. Then either $G_L$ is {2-}transitive on $\PP \setminus L$, or $G_L\leq G_\s$ and $|L|(|L|-1)$ is a divisor of $|\PP|-1$. Moreover, if $L=\{\a, \b\}$, then either $|\PP|=4$ and $G=A_4$, or $G_{\a, \b}$ is transitive on $\PP \setminus\{\a, \b\}$, or $|\PP|$ is odd and $G_{\a, \b}\leq G_{\s}$ is $\frac{1}{2}$-transitive on $\PP \setminus\{\a, \b, \s\}$.
\end{lemma}


\proof
Since $\Om$ is feasible, $G_{\s, L}=G_{\s, \PP \setminus L}$ is transitive on $(\PP \setminus L) \setminus \{\s\}$. If $G_L\nleq G_\s$, then $G_L$ is 2-transitive on $\PP \setminus L$. Assume that $G_L\leq G_\s$ in the sequel, so that $|G_L|$ is a divisor of $|G_\s|$. There are $b_\DD=|\PP|(|\PP|-1)/(|L|(|L|-1))$ blocks in $\DD$. Since $|G_L|=|G|/b_\DD$ and $|G_\s|=|G|/|\PP|$, $|\PP|$ {needs} to be a divisor of $b_\DD$ and hence $|L|(|L|-1)$ is a divisor of $|\PP|-1$.

Suppose that $L=\{\a, \b\}$ has order $2$ (that is, $\DD$ is improper). First assume $G_{\{\a, \b\}}\nleq G_\s$, and thus $G_{\{\a, \b\}}$ is $2$-transitive on $\PP \setminus\{\a, \b\}$. If $G_{\{\a, \b\}}$ is not faithful on $\PP \setminus\{\a, \b\}$, then the transposition $(\a\,\b)$ is contained in $G$ and thus $G$ contains all transpositions by the 2-transitivity of $G$ on $\PP$, which implies $G=\Sym(\PP)$ and hence $G_{\a, \b}$ is transitive on $\PP \setminus\{\a, \b\}$. We may assume $G_{\{\a, \b\}}$ is faithful on $\PP \setminus\{\a, \b\}$. If $G_{\a, \b}\neq1$, then $G_{\a, \b}$ is transitive on $\PP \setminus\{\a, \b\}$ by \cite[Theorem 9.9]{Wielandt} and \cite[\textsection 39]{Manning1921} as $G_{\a, \b}$ is a normal subgroup of $G_{\{\a, \b\}}$. If $G_{\a, \b}=1$, then $|G_{\{\a, \b\}}|=2$ and $|\PP|=4$ as $G_{\{\a, \b\}}$ is $2$-transitive on $\PP \setminus\{\a, \b\}$, which implies $G=A_4$ by the $2$-transitivity of $G$ on $\PP$. Next assume $G_{\{\a, \b\}}\leq G_\s$. Then $|\PP|$ is odd and $G_{\a, \b}\unlhd G_{\{\a, \b\}}=G_{\{\a, \b\}, \s}$ is $\frac{1}{2}$-transitive on $\PP \setminus\{\a, \b, \s\}$.
\qed

\section{Proof of Theorem \ref{thm:main theorem}: Part (a)}
\label{sec:parta}

By Lemma \ref{lem:two design}, to prove Theorem \ref{thm:main theorem} it suffices to classify all $G$-flag graphs of the complements of $(G, 2)$-point-transitive linear spaces $\DD$ with $u = |\BB|$ points. We will achieve this goal in this and the next sections for proper and improper $\DD$, respectively.

\subsection{Doubly point-transitive proper linear spaces}
\label{sec:doubly point-transitive linear spaces}

All proper $(G, 2)$-point-transitive linear spaces $\DD$ are known \cite[Theorem 1]{Kantor85}. In the case when $G$ is almost simple, $(\DD, G)$ is one of the following:
\begin{itemize}[topsep=0.5ex]
\setlength\itemsep{-0.3em}
\item[(i)] $\DD=\PG(3,2)$, $G = A_7$;
\item[(ii)] $\DD$ is the $2$-$(q^3+1, q+1, 1)$ design $U_{H}(q)$ (Hermitian unital) associated with $\PSU(3, q)$, and $\soc(G) = \PSU(3, q)$, where $q > 2$ is a prime power;
\item[(iii)] $\DD$ is the $2$-$(q^3+1, q+1, 1)$ design $U_{R}(q)$ (Ree unital \cite{Lu}) associated with the Ree group $\R(q)$, and $\soc(G) =\R(q)'$, where $q = 3^{2e+1} \geq 3$;
\item[(iv)] $\DD = \PG(d-1,q)$, $\PSL(d, q) \le G \le \PGammaL(d, q)$ ($\soc(G) = \PSL(d, q)$), $d \ge 3$, $q=p^n$ for some prime $p$.
\end{itemize}

In each case above, except for $G =\R(3)$ ($\cong \PGammaL(2, 8)$), $\soc(G)$ is also $2$-transitive on the set of points of $\DD$.

In the case when $G$ is affine, $\soc(G)$ is of order $u = p^d$, where $d \ge 1$ and $p$ is a prime, and we may identify $G$ with a group of affine transformations $\bv \mapsto \bv^g + \bc$ of $\mathbb{F}_p^d$, where $\bv, \bc \in \mathbb{F}_p^d$ and $g$ is in the stabilizer $G_{\bf 0} \le \GL(d, p)$ of ${\bf 0} \in \mathbb{F}_p^d$ in $G$. The proof of \cite[Theorem 1]{Kantor85} implies that all possibilities for $(G, \DD)$ are as follows, where $q$ is a power of $p$.

\begin{itemize}[topsep=0.5ex]
\setlength\itemsep{-0.3em}
\item[(v)] $G \le \AGammaL(1, u)$, $\DD$ is an affine space;
\item[(vi)] $\SL(n, q) \unlhd G_{\bf 0}$, $u = q^n$, $n\geq2$;
\item[(vii)] $\Sp(n, q) \unlhd G_{\bf 0}$, $u = q^n$, $n \ge 4$ is even;
\item[(viii)] $G_2(q)' \unlhd G_{\bf 0}$, $u = q^6$, $q$ is even;
\item[(ix)] $\SL(2,3) \unlhd G_{\bf 0}$ or $\SL(2,5) \unlhd G_{\bf 0}$, $\DD = \AG(2,p)$, $u=p^2$, $p = 5, 7, 11, 19, 23, 29$ or $59$;
\item[(x)] $G_{\bf 0}$ contains a normal extraspecial subgroup $E$ of order $2^5$, $\DD = \AG(4,3)$, $u=3^4$;
\item[(xi)] $G_{\bf 0}$ contains a normal extraspecial subgroup $E$ of order $2^5$, $\DD$ is the unique ``exceptional nearfield affine plane'' of order $9$ having $3^4$ points and $3^2 \cdot (3^2 + 1)$ lines, with $3^2$ points on each line and $3^2 + 1$ lines through each point (see \cite{Foulser} and \cite[pp.33-34, 229-232]{Dembowski}), and $u=3^4$;
\item[(xii)] $\SL(2, 5) \unlhd G_{\bf 0}$, $\DD$ is $\AG(2,9)$, the ``exceptional nearfield affine plane'' as in (xi), or $\AG(4,3)$ as in (x), $u=3^4$;
\item[(xiii)] $G_{\bf 0} = \SL(2, 13)$, $\DD = \AG(6,3)$,  $u=3^6$;
\item[(xiv)] $G_{\bf 0} = \SL(2, 13)$, $\DD$ is the Hering affine plane of order $27$ having $3^6$ points and $3^3 \cdot (3^3 + 1)$ lines, with $3^3$ points on each line and $3^3 + 1$ lines through each point (see \cite{Hering27} and \cite[p.236]{Dembowski}), and $u=3^6$;
\item[(xv)] $G_{\bf 0} = \SL(2, 13)$, $\DD$ is one of the two Hering designs \cite{Hering} of order $90$ having $3^6$ points and $81 \cdot 91$ lines, with $9$ points on each line and $91$ lines through each point, and $u=3^6$.
\end{itemize}

\subsection{Almost simple case}
\label{sec:almost simple}

We first prove that among cases (i)--(iv) above a feasible $G$-orbit on the flag set of $\overline{\DD}$ exists only in case (iv) with $d=3$. We also classify the corresponding $G$-flag graphs and establish their properties.

Consider case (i) first. Suppose that $\Om = (\s, \PP \setminus L)^G$ is a feasible $G$-orbit on the set of flags of $\overline{\DD}$, where $(\s, L)$ is an antiflag of {the $2$-$(15, 3, 1)$ design $\DD = \PG(3,2)$}. If $G_L$ is $2$-transitive on $\PP \setminus L$, then $|\PP \setminus L|(|\PP \setminus L|-1)=12\cdot11$ divides $G_L$, a contradiction. {Thus, by Lemma \ref{lem:two 2-transitive actions},} $|L|(|L|-1)=6$ divides $|\PP|-1=14$, which is a contradiction. Hence there is no feasible $G$-orbit on the set of flags of $\overline{\DD}$.

Consider case (ii) next. Suppose that $\Om = (\s, \PP \setminus L)^G$ is a feasible $G$-orbit on the set of flags of $\overline{\DD}$, where $(\s, L)$ is an antiflag of $\DD$. Since $|L|(|L|-1)=(q+1)q$ cannot divide $|\PP|-1=q^3$, by Lemma \ref{lem:two 2-transitive actions}, $G_L$ is $2$-transitive on $\PP \setminus L$ and thus $|\PP \setminus L|(|\PP \setminus L|-1)=(q^3-q)(q^3-q-1)$ divides $G_L$. We know {that} $\PSU(3,q)$ is of index $\gcd(3, q+1)$ in $\PGU(3, q)$ (\cite[p.249]{Dixon-Mortimer}), $|\GU(3, q)|=q^3(q+1)(q^2-1)(q^3+1)$ and $|\PGU(3, q)|=|\GU(3, q)|/(q+1)$ (see \cite[p.66]{Wilson}). Moreover, there are $b_\DD=(q^3+1)q^3/((q+1)q)$ blocks in $\DD$, and $|G_L|=|G|/{b_\DD}$ divides $d\cdot|\PGU(3, q)|/{b_\DD}=dq(q+1)(q^2-1)$ as $[G:\PSU(3,q)]$ is a divisor of $d\cdot\gcd(3, q+1)$ (see \cite[p.197, Table 7.4]{Cameron99}). But $(q^3-q)(q^3-q-1)>dq(q+1)(q^2-1)$ when $q>2$, which implies that $G_L$ cannot be $2$-transitive on $\PP \setminus L$. Hence there is no feasible $G$-orbit on the set of flags of $\overline{\DD}$.

In case (iii), suppose that $\Om = (\s, \PP \setminus L)^G$ is a feasible $G$-orbit on the set of flags of $\overline{\DD}$, where $(\s, L)$ is an antiflag of $\DD$. Since $|L|(|L|-1)=(q+1)q$ cannot divide $|\PP|-1=q^3$, by Lemma \ref{lem:two 2-transitive actions}, $G_L$ is $2$-transitive on $\PP \setminus L$ and thus $|\PP \setminus L|(|\PP \setminus L|-1)=(q^3-q)(q^3-q-1)$ divides $G_L$. There are $b_\DD=(q^3+1)q^3/((q+1)q)$ blocks in $\DD$, and thus $|G_L|=|G|/{b_\DD}$ divides $(2e+1)|\R(q)|/{b_\DD}=(2e+1)q(q^2-1)$ as $[G:\R(q)]$ is a divisor of $2e+1$ (see \cite[p.197, Table 7.4]{Cameron99}). But $(q^3-q)(q^3-q-1)>(2e+1)q(q^2-1)$ when $q>2$, which implies that $G_L$ cannot be $2$-transitive on $\PP \setminus L$. Therefore, there is no feasible $G$-orbit on the set of flags of $\overline{\DD}$.

The rest of this subsection is devoted to case (iv). Suppose that $\Om = (\s, \PP \setminus L)^G$ is a feasible $G$-orbit on the set of flags of $\overline{\DD}$, where $(\s, L)$ is an antiflag of $\DD$.
It is known (see \cite[1.4.24]{Dembowski}) that for any integer $s$ between $1$ and $d-1$, $G$ is transitive on the set of ordered $(s+1)$-tuples of independent points of $\PG(d-1, q)$, where $s+1$ points of $\PG(d-1, q)$ are said to be independent if they do not lie on any $(s-1)$-flat of $\PG(d-1, q)$. In particular, fixing a line $L$ of $\PG(d-1, q)$, say, the line $L = \la \bfe_1, \bfe_2 \ra$ corresponding to the plane of $\mathbb{F}_q^d$ spanned by  $\bfe_1$ and $\bfe_2$, the stabilizer $G_{\a_1, \a_2}$ where $\a_1 = \la \bfe_1 \ra, \a_2 = \la \bfe_2 \ra$, is transitive on $\PP \setminus L$. Since $G_{\a_1, \a_2} \le G_L$, $G_L$ is also transitive on $\PP \setminus L$ and thus $\Om := (\s, \PP \setminus L)^G$ is the set of flags of $\overline{\DD}$, where $\s = \la \bfe_d \ra$.

A typical element of $G\leq \PGammaL(d, q)$ is
$$
t(A, j):\la\bx\ra\mapsto \la\bx^{\zeta^j} A\ra,
$$
where $A\in \GL(d, q)$ (we always identify scale multiples in the following discussion), $\bx\in \mathbb{F}_q^d$, $j$ is an integer and $\zeta:\mathbb{F}_q\rightarrow\mathbb{F}_q$, $z\mapsto z^p$ is the Frobenius map {acting componentwise on $\bx$}. It can be verified that if $t(A, j)\in G_{\s, L}$ then $A$ is of the form
$$
\bmat{A_{11} & 0 & 0\\A_{21} & A_{22} & \bb^T\\ 0 & 0 & 1},
$$
where $A_{11} \in \GL(2, q), A_{22} \in \GL(d-3, q)$, $A_{21}$ is a $(d-3) \times 2$ matrix and $\bb^T$ is a $(d-3) \times 1$ matrix over $\mathbb{F}_q$. (Note that $A_{22}$ vanishes when $d=3$.) Since any element of $G_{\s, L}$ maps $\bfe_1 + \bfe_d$ to a vector of $\mathbb{F}_q^d$ in the 3-dimensional subspace $\la \bfe_1, \bfe_2, \bfe_d \ra$, if $d \ge 4$, then $G_{\s, L}$ cannot be transitive on $(\PP \setminus L) \setminus \{\s\}$ and $\Om$ is not feasible. On the other hand, if $d = 3$, then $G_{\s, L}$ is transitive on $(\PP \setminus L) \setminus \{\s\}$ and so $\Om$ is feasible.

Assume $d = 3$ in the sequel. As $\overline{\DD}$ is a $2$-$(v_{\overline{\DD}}, k_{\overline{\DD}}, \lambda_{\overline{\DD}})$ design with $b_{\overline{\DD}}$ blocks, we obtain $v_{\overline{\DD}}=v_\DD=q^2+q+1$, $k_{\overline{\DD}}=v_{\DD}-k_{\DD}=q^2$, $b_{\overline{\DD}}=b_\DD=q^2+q+1$ and $\l_{\overline{\DD}}=q(q-1)$. Thus $|\Om|=b_{\overline{\DD}}k_{\overline{\DD}} = q^2 (q^2+q+1)$ and $|\Om(\s)|=|\Om|/{v_{\overline{\DD}}}=q^2$. We claim that every $G$-orbit on $\C(\overline{\DD}, \Om)$ is self-paired. In fact, for any $G$-orbit $\Psi:=((\eta, \PP \setminus L), (\xi, \PP \setminus N))^G$ on $\C(\overline{\DD}, \Om)$, we {may} assume $\eta=\la \bfe_3\ra$ and $L=\la \bfe_1, \bfe_2\ra$. Since $G_{\eta, L}$ is transitive on $(\PP \setminus L) \setminus \{\eta\}$, we {may} further assume $\xi=\la \bfe_1+\bfe_3\ra$. Thus
$$
\Psi = ((\la \bfe_3\ra, \PP \setminus \la \bfe_1, \bfe_2\ra), (\la \bfe_1+\bfe_3\ra, \PP \setminus N))^G.
$$
We know {that}
\begin{align*}
N=\la (0, 1, h), (1, 0, f)\ra \text{\;for some\;}  h, f \in \mathbb{F}_q \text{\;with\;} f\neq1.
\end{align*}
Set $g:=t(C, 0)$, where
$$
C:=\bmat{-1 & 0 & -f\\\frac{-h}{f-1} & \frac{1}{f-1} & -h\\ 1 & 0 & 1}.
$$
Then $g\in \PSL(3,q)\leq G$, $(\eta, \xi)^g=(\xi, \eta)$ and $(L, N)^g=(N, L)$. Therefore, $\Psi$ is self-paired as claimed and defines a $G$-flag graph $\Ga(\overline{\DD}, \Om, \Psi)$.

We now prove that in $\Ga(\overline{\DD}, \Om, \Psi)$ there are exactly $q(q-1)$ vertices in $\Om(\eta)$ that have neighbours in $\Om(\xi)$, and they are the vertices $(\eta, \PP \setminus M) \in \Om(\eta)$, where $\PP \setminus M$ is running over all $q(q-1)$ blocks of $\overline{\DD}$ containing both $\eta$ and $\xi$. In fact, if $(\eta, \PP \setminus M)$ is a vertex in $\Om(\eta)$ which has a neighbour in $\Om(\xi)$, then by \eqref{equ:C(D,Om)}, $\PP \setminus M$ contains both $\eta$ and $\xi$. Conversely, let $\PP \setminus M$ be a block of $\overline{\DD}$ containing both $\eta$ and $\xi$. Then there exist $x, y \in \mathbb{F}_q$ with $y\neq1$ such that $M = \la (0, 1, x), (1, 0, y)\ra$. Set $z:=t(R, 0)$, where
$$
R:=\bmat{\frac{1}{1-y} & 0 & \frac{y}{1-y}\\x & 1-y & x\\ 0 & 0 & 1}.
$$
Then $z\in G_{\eta, \xi}$ and $L^z=M$. Since $(\eta, \PP \setminus L)$ is adjacent to $(\xi, \PP \setminus N)$, it follows that $(\eta, \PP \setminus M)$ is adjacent to $(\xi, \PP \setminus N^z)\in \Om(\xi)$.

We now determine the number of neighbours of $(\eta, \PP \setminus L)$ in $\Om(\xi)$ in the $G$-flag graph $\Ga(\overline{\DD}, \Om, \Psi)$. This number, denoted by $\ell$, is also the number of neighbours in $\Om(\xi)$ of each of the $q(q-1)$ vertices in $\Om(\eta)$ that are adjacent to at least one vertex in $\Om(\xi)$. Set $J:=G_{\eta, \xi, L}$. {Then $(\xi, \PP \setminus N)^{J}$} is the set of vertices in $\Om(\xi)$ adjacent to $(\eta, \PP \setminus L)$ in $\Ga(\overline{\DD}, \Om, \Psi)$, and $\ell = |(\xi, \PP \setminus N)^{J}| = |J/J_N|$, where $J_N = \{t(A, j)\in J: N^{t(A, j)}=N\}$ is the stabilizer of $N$ in $J$.

In order to determine $\ell$ we introduce some notions first. Let $\widehat{c}: \mathbb{F}_q \rightarrow \mathbb{F}_q, x\mapsto cx$ be the scalar multiplication by $c\in\mathbb{F}_q^\times$. Choose $\omega$ to be a fixed primitive element of $\mathbb{F}_q$. Then $\widehat{\omega}$ generates the multiplicative group $\GL(1, q)$ and $\GammaL(1, q)=\la \widehat{\om}, \zeta\ra$. Foulser \cite{Foulser64} gives a standard generating set for each subgroup of $\GammaL(1, q)$.

\begin{definition}
\label{def:standard form}
{\em (\cite[Definition 4.5]{LiLimPraeger})} {\em Every subgroup $M$ of $\GammaL(1,p^{n})$ can be uniquely presented in the form $M=\la \widehat{\om}^t, \zeta^s \widehat{\om}^e\ra$ such that the following conditions all hold:
\begin{itemize}
\item[\rm (F1)] $t>0$ and $t\mid (p^{n}-1)$;
\item[\rm (F2)] $s>0$ and $s\mid n$;
\item[\rm (F3)] $0\leq e<t$ and $t\mid e(p^{n}-1)/(p^s-1)$.
\end{itemize}
The presentation $M=\la \widehat{\om}^t, \zeta^s \widehat{\om}^e\ra$ satisfying {(F1)}--{(F3)} is said to be in {\em standard form} with {\em standard parameters} $(t, e, s)$.}
\end{definition}

The reader is referred to \cite[Lemma 4.1]{Foulser64} for the existence and uniqueness of a standard form for each subgroup of $\GammaL(1, p^{n})$, and to \cite{FFXZ, FFXZ2, Foulser64} for information about the structure of subgroups of $\GammaL(1,p^{n})$. Set
$$
H:=\{t(D_c, j): t(D_c, j)\in G\}\text{\; and \;} E:=\{t(T_a, 0): t(T_a, 0)\in G\},
$$
where
$$
D_c:=\bmat{1 & 0 & 0\\0 & c & 0\\ 0 & 0 & 1} \text{\; and \;} T_a:=\bmat{1 & 0 & 0\\a & 1 & 0\\ 0 & 0 & 1}.
$$
It can be shown that $G=\PSL(3, q)\rtimes H$, and for any $t(A, j)\in J$ we have $A=T_a D_c$ for some $a, c \in\mathbb{F}_q$ with $c\neq0$ which implies $J=E\rtimes H$ (if $t(A, j)\in J_N$, then in addition $f^{\zeta^j}=f$ and $h^{\zeta^j}=af+ch$). Define
$$
\Lambda: H\rightarrow \GammaL(1,q),\;\, t(D_c, j)\mapsto t(c,j),
$$
{where $t(c,j): y \mapsto cy^{\zeta^j}$ for $y \in \mathbb{F}_q$}. Then $\Lambda$ is a homomorphism and is injective. Thus $H$ is isomorphic to a subgroup of $\GammaL(1,q)$.

Suppose that $\Lambda(H)=\la \widehat{\om}^t, \zeta^s \widehat{\om}^e\ra$ is in standard form. Then $t(A,j)\in J$ implies $s\mid j$, and $|J|=q(q-1)n/(st)$. Moreover, since $t(D_{c^3}, 0)\in \PSL(3,q)\cap J$ for $c \in\mathbb{F}_q^\times$, we obtain $\la\widehat{\omega}^3\ra\leq \Lambda(H)\cap \GL(1, q)$ and thus $t\mid 3$ ($t=1$ if and only if $\PGL(3,q)\leq G$). Let
$$
\rho:=\min\{i>0: f^{p^i}=f\}.
$$
If $t(A, j)\in J_N$, then we can see that $\rho\mid j$. Hence $t(A, j)\in J_N$ implies $\lcm(\rho, s)\mid j$. Set  
$$
\mu:=\min\{i>0: n \text{\;divides\;} i\cdot\lcm(\rho, s)\}.
$$
Since $s\mid n$ and $\rho\mid n$ by definition, we have $\lcm(\rho, s)\mid n$ and thus $\mu=n/\lcm(\rho, s)$. Therefore, it can be verified that $|J_N|=|J|$ when $f=h=0$ and $|J_N|=\mu(q-1)/t$ when $f\neq0$. In the case when $f=0$ and $h\neq0$, assuming $h^{-1}=\omega^{\theta}$, where $0<\theta<q$, let $m$ be the smallest positive integer such that $t\mid (e+\theta (p^s-1))\cdot \frac{p^{sm}-1}{p^s-1}$. Then we can verify that $m$ divides $n/s$, $\Lambda(H)_{h^{-1}}$ is cyclic of order $n/(ms)$, and $|J_N|=qn/(ms)$.

\medskip
(A1) $f=h=0$. In this case, $N = \la e_1, e_2\ra$ and $\ell = |J|/|J_N|=1$. The $G$-flag graph $\Ga$ has order $q^2(q^2+q+1)$ and valency $q^2-1$, and $\Ga[\Om(\eta), \Om(\xi)]\cong q(q-1)\cdot K_2$.  

\medskip
(A2) $f=0$ and $h\neq0$. In this case, $\ell = |J|/|J_N| = (q-1)m/t=q-1$ if $\PGL(3,q)\leq G$ (in this case $t=1$ and thus $m=1$) or $(q-1)m/3$ if $\PGL(3,q)\nleq G$. The $G$-flag graph $\Ga$ has $q^2(q^2+q+1)$ vertices, the number of vertices in $\Om(\eta)$ that have neighbours in $\Om(\xi)$ is $q(q-1)$, and the valency of $\Ga$ is $(q^2-1)(q-1)$ if $\PGL(3,q)\leq G$ or $(q^2-1)(q-1)m/3$ if $\PGL(3,q)\nleq G$.

\medskip
(A3) $m\neq0$ and $m\neq1$. In this case, $\ell = |J|/|J_N| = q\cdot\lcm(\rho, s)/s$. The $G$-flag graph $\Ga$ has $q^2(q^2+q+1)$ vertices, the number of vertices in $\Om(\eta)$ that have neighbours in $\Om(\xi)$ is $q(q-1)$, and the valency of $\Ga$ is $(q^2-1)q\cdot\lcm(\rho, s)/s$.

\subsection{Affine case}
\label{sec:affine}

We now prove that in cases (v)--(xv) in Section \ref{sec:doubly point-transitive linear spaces} there is no feasible $G$-orbit on the set of flags of $\overline{\DD}$. Suppose otherwise. Let $\Om = (\s, \PP \setminus L)^G$ be such a $G$-orbit, where $(\s, L)$ is an antiflag of $\DD$.

In cases (ix)--(xv), since $|L|(|L|-1)$ {does not} divide $|\PP|-1$, by Lemma \ref{lem:two 2-transitive actions}, $G_L$ is $2$-transitive on $\PP \setminus L$. However, we can check that $(|\PP|-|L|)(|\PP|-|L|-1)$ cannot divide $|G_L|=|L|\cdot|G_{\bzero, L}|$, which is a contradiction.

In case (v), since $|L|(|L|-1)$ cannot divide $|\PP|-1$, by Lemma \ref{lem:two 2-transitive actions}, $G_L$ is $2$-transitive on $\PP \setminus L$. Suppose that $|L|=q$ and $p^d=q^n$. Then $b_\DD=\frac{q^n(q^n-1)}{q(q-1)}$ and $|G_L|=|G|/{b_\DD}$. Since $G\leq\AGammaL(1, q^n)$, $|G|$ is a divisor of $|\AGammaL(1, q^n)|=q^n(q^n-1)d$. Hence $(q^n-q)(q^n-q-1)$ divides ${q^n(q^n-1)d}/{b_\DD}$ and $dq(q-1)/((q^n-q)(q^n-q-1))$ is an integer, which implies $d\geq (q^n-q)/(q-1)$ as $\DD$ is proper (thus $q>2$ and $n>1$) by our assumption. Therefore, $d>(q^n-q)/q=q^{n-1}-1$ and $d+1>q^{n-1}=p^{d(1-1/n)}\geq p^{d/2}$. We obtain $(p, n, d)=(2, 2, 4)$. But $dq(q-1)/((q^n-q)(q^n-q-1))=48/(12\cdot11)$ is not an integer, a contradiction.

In cases (vi), (vii) and (viii), $\DD$ contains a block $L\subseteq\la\bx\ra$ and $L$ is a vector space over $\mathbb{F}_p$, where $\bzero\neq\bx\in L$ \cite[Section 4]{Kantor85}. Without loss of generality, {we may} assume that $\bx=\bfe_1$.  Since $G\leq \AGammaL(n,q)$, a typical element of $G$ is $\t(A, \bc, \vp):\bz\mapsto\bz^\vp A+\bc$, where $A\in \GL(n, q)$, $\bc\in\mathbb{F}_q^n$ and $\vp$ is a field automorphism of $\mathbb{F}_q$ acting componentwise on $\bz$. Let $\s:=\bfe_n$. It can be verified that if $\t(A, \bc, \vp)\in G_{\s, L}$ then $\bc\in \la\bfe_1\ra$ {and} $A$ is of the form
$$
\bmat{a_{11} & \bzero & 0\\A_{21} & A_{22} & A_{23}\\ -c & \bzero & 1},
$$
where $a_{11}\neq0$, $ A_{22} \in \GL(n-2, q)$, $A_{21}$ and $A_{23}$ are both $(n-2) \times 1$ matrix over $\mathbb{F}_q${,} and $\bc=(c,0,\ldots,0)${.} (Note that $A_{22}$ vanishes when $n=2$.) Since {every} element of $G_{\s, L}$ maps $\bfe_1+\bfe_n$ to a vector of $\mathbb{F}_q^n$ in the $2$-dimensional subspace $\la\bfe_1, \bfe_n\ra$, if $n>2$, then $G_{\s, L}$ cannot be transitive on $(\PP \setminus L) \setminus \{\s\}$. In the case when $n=2$, if $q>2$, then $(1,1)$ cannot be mapped to $(0, z)$ by {any element} of $G_{\s, L}$, where $z\neq1$. Hence $|\PP|=|\mathbb{F}_2^2|=4$ and $\DD$ is improper, which is a contradiction.

\medskip
So far we have completed the proof of part (a) of Theorem \ref{thm:main theorem}.

\section{Proof of Theorem \ref{thm:main theorem}: Part (b)}
\label{sec:trivial}

In this section we assume that $\DD$ is an improper $(G, 2)$-point-transitive linear space on $\PP$ with $u := |\PP|$ points. Then $\overline{\DD}$ is a $2$-$(u, u-2, (u-2)(u-3)/2)$ design. We may identify $\DD$ with the $G$-arc transitive complete graph with vertex set $\PP$, and identify $\overline{\DD}$ with the complete $(u-2)$-uniform hypergraph with vertex set $\PP$ and hyperedge set ${\PP \choose u-2}$. Let $\a, \b, \s \in \PP$ be distinct points so that $(\s, \PP \setminus \{\a, \b\})$ is a flag of $\overline{\DD}$. Suppose that $\Om = (\s, \PP \setminus \{\a, \b\})^G$ is feasible and $\Psi=((\s, \PP \setminus \{\a, \b\}), (\t, \PP \setminus \{\g, \d\}))^G$ is a self-paired $G$-orbit on $\C(\overline{\DD}, \Om)$. Then $\Om$ and $\Psi$ produce a $G$-flag graph $\Ga:=\Ga(\overline{\DD}, \Om, \Psi)$. By Lemma \ref{lem:properties about the G-flag graph}, the number of vertices of $\Ga$ is $u\cdot|\Om(\s)|=u\cdot|\{\a,\b\}^{G_{\s}}|$, the number of vertices in $\Om(\s)$ which have neighbours in $\Om(\t)$ is $|\{\a,\b\}^{G_{\s,\t}}|$, each vertex in $\Om(\s)$ is adjacent to $|\{\g,\d\}^{H}|$ vertices in $\Om(\t)$, and the valency of $\Ga$ is $(u-3)\cdot|\{\g,\d\}^{H}|$, where $H=G_{\s,\t,\{\a,\b\}}$.

\subsection{Almost simple case}
\label{sec:almost simple 2-transitive groups}

Suppose that $G\leq\Sym(\PP)$ is $2$-transitive on $\PP$ of degree $u=|\PP|$ with $\soc(G)$ a nonabelian simple group. Then $\soc(G)$ and $u$ are as follows (\cite{Kantor85}, \cite[p.196]{Cameron99}, \cite{Cameron81}):
\begin{itemize}[topsep=0.5ex]
\setlength\itemsep{-0.3em}
\item[(i)] $\soc(G) = A_u$, $u \geq 5$;
\item[(ii)] $\soc(G) = \PSL(d, q)$, $d \geq 2$, $q$ is a prime power and $u = (q^d-1)/(q-1)$, where $(d, q)\neq$ $(2,2), (2,3)$;
\item[(iii)] $\soc(G) = \PSU(3, q)$, $q \geq 3$ is a prime power and $u = q^3 + 1$;
\item[(iv)] $\soc(G) = \Sz(q)$, $q = 2^{2e+1} > 2$ and $u = q^2+1$;
\item[(v)] $\soc(G) = \R(q)'$, $q = 3^{2e+1} $ and $u = q^3+1$;
\item[(vi)] $G = \Sp_{2d}(2)$, $d \geq 3$ and $u = 2^{2d-1} \pm 2^{d-1}$;
\item[(vii)] $G=\PSL(2,11)$, $u = 11$;
\item[(viii)] $G = \HS$, $u = 176$;
\item[(ix)] $G = \Co_3$, $u = 276$;
\item[(x)] $G = M_{11}$, $u = 12$;
\item[(xi)] $G=A_7$, $u=15$;
\item[(xii)] $\soc(G) = M_u$, $u = 11,12,22,23,24$.
\end{itemize}

\delete
{
\begin{description}
\item[\rm (i)] $\soc(G) = A_u$, $u \geq 5$;

\item[\rm (ii)] $\soc(G) = \PSL(d, q)$, $d \geq 2$, $q$ is a prime power and $u = (q^d-1)/(q-1)$, where $(d, q)\neq$ $(2,2), (2,3)$;

\item[\rm (iii)] $\soc(G) = \PSU(3, q)$, $q \geq 3$ is a prime power and $u = q^3 + 1$;

\item[\rm (iv)]  $\soc(G) = \Sz(q)$, $q = 2^{2e+1} > 2$ and $u = q^2+1$;

\item[\rm (v)] $\soc(G) = \R(q)'$, $q = 3^{2e+1} $ and $u = q^3+1$;

\item[\rm (vi)] $G = \Sp_{2d}(2)$, $d \geq 3$ and $u = 2^{2d-1} \pm 2^{d-1}$;

\item[\rm (vii)] $G=\PSL(2,11)$, $u = 11$;

\item[\rm (viii)] $G = \HS$, $u = 176$;

\item[\rm (ix)] $G = \Co_3$, $u = 276$;

\item[\rm (x)] $G = M_{11}$, $u = 12$;

\item[\rm (xi)] $G=A_7$, $u=15$;

\item[\rm (xii)] $\soc(G) = M_u$, $u = 11,12,22,23,24$.
\end{description}
}

In case (i), $G_{\{\a, \b\}, \s}$ is transitive on $\PP \setminus\{\a,\b, \s\}$ and thus $\Om = (\s, \PP \setminus \{\a, \b\})^G$ is feasible. Let $\Psi=((\s, \PP \setminus \{\a, \b\}), (\t, \PP \setminus \{\g, \d\}))^G$ be a $G$-orbit on $\C(\overline{\DD}, \Om)$. Then $\Psi$ is self-paired as $A_u\leq G$, and thus $\Om$ and $\Psi$ produce a $G$-flag graph $\Ga:=\Ga(\overline{\DD}, \Om, \Psi)$. We can see that $|\Om(\s)|=(u-1)(u-2)/2$, $|\Om|=u(u-1)(u-2)/2$, and the number of vertices in $\Om(\s)$ which have neighbours in $\Om(\t)$ is $(u-2)(u-3)/2$. Set $\ell:=|\{\a, \b\}\cap\{\g, \d\}|$. If $\ell=2$, then $(\s, \PP \setminus \{\a, \b\})$ is adjacent in $\Ga$ to only one vertex $(\t, \PP \setminus \{\a, \b\})$ in $\Om(\t)$. If $\ell=1$, then $(\s, \PP \setminus \{\a, \b\})$ is adjacent in $\Ga$ to one vertex in $\Om(\t)$ when $G=A_5$, two vertices in $\Om(\t)$ when $G=S_5$, two vertices in $\Om(\t)$ when $G=A_6$, four vertices in $\Om(\t)$ when $G=S_6$ and $2u-8$ vertices in $\Om(\t)$ when $u>6$. The case $\ell=0$ occurs only when $u>5$, and in this case $(\s, \PP \setminus \{\a, \b\})$ is adjacent to $(u-4)(u-5)/2$ vertices in $\Om(\t)$.

In case (ii) with $d\geq3$, $G_{\a,\b}$ has orbit-lengths $q-1$ and $u-(q+1)$ on $\PP \setminus\{\a,\b\}$, and thus it is not transitive on $\PP \setminus\{\a,\b\}$. By Lemma \ref{lem:two 2-transitive actions}, this case does not produce any $G$-flag graph.

In case (ii) with $d=2$ and $q=p^e\geq5$ (where $p$ is a prime), let $\a=\la \bfe_1\ra$ and $\b=\la \bfe_2\ra$, where $\bfe_1=(1,0)$ and $\bfe_2=(0,1)$. Suppose that there exists a feasible $G$-orbit on the flag set of $\overline{\DD}$. Since each orbit of $G_{\a,\b}$ on $\PP \setminus\{\a,\b\}$ has length at least $(q-1)/2\geq2$ (see \cite{Kantor85}), $G_{\{\a,\b\}}$ must be 2-transitive on $\PP \setminus\{\a,\b\}$ by Lemma \ref{lem:two 2-transitive actions}. Hence $(q-1)(q-2)$ divides $J_{\{\a,\b\}}$, where $J:=\PGammaL(2,q)$, as $G\leq\PGammaL(2,q)$, and we obtain that $q-2$ divides $2e$. This can happen only when $p=2$ and $e=3$. It can be verified that $J_{\{\a,\b\}, \s}$ is transitive on $\PP \setminus\{\a,\b, \s\}$, where $\s:=\la(1,1)\ra$, and $\Om:=(\s, \PP \setminus\{\a,\b\})^J$ is feasible and is the set of all flags of $\overline{\DD}$. Let $\Psi$ be a $G$-orbit on $\C(\overline{\DD}, \Om)$. Since $G_{\{\a, \b\}}$ is 2-transitive on $\PP \setminus\{\a, \b\}$, we may assume
$$\Psi= ((\la\bfe_1\ra, \PP\setminus\{\la (1, 1)\ra, \la(1, u)\ra\}), (\la\bfe_2\ra, \PP\setminus \{\la(1, x)\ra, \la(1, y)\ra\}))^G,$$
where $u$ is a generator of $\mathbb{F}_8^\times$, $x, y\neq0$ and $x\neq y$. To prove $\Psi$ is self-paired, we need to find some $g\in J$ such that $g$ interchanges $(\la\bfe_1\ra, \PP\setminus\{\la (1, 1)\ra, \la(1, u)\ra\})$ and $(\la\bfe_2\ra, \PP\setminus \{\la(1, x)\ra, \la(1, y)\ra\})$. In fact, if $x=u^{\vp}y$, where $\vp\in \Aut(\mathbb{F}_8)$, then we set $g:\PP\rightarrow\PP$, $\la\bc\ra\mapsto \la\bc^\vp A\ra$, where $A=\bmat{0 & x\\ 1 & 0}$; if $y=u^{\vp}x$, where $\vp\in \Aut(\mathbb{F}_8)$, then we set $g:\PP\rightarrow\PP$, $\la\bc\ra\mapsto \la\bc^\vp A\ra$, where $A=\bmat{0 & y\\ 1 & 0}$. For each self-paired $G$-orbit $\Psi$ on $\C(\overline{\DD}, \Om)$, the corresponding $G$-flag graph $\Ga:=\Ga(\overline{\DD}, \Om, \Psi)$ has $|\Om|=252$ vertices, the number of vertices in $\Om(\la\bfe_1\ra)$ which have neighbours in $\Om(\la\bfe_2\ra)$ is $21$, and $\Ga[\Om(\la\bfe_1\ra), \Om(\la\bfe_2\ra)]\cong 21\cdot K_2$.


In cases (iii)--(iv), case (v) with $q>3$ and cases (vi)--(ix), $G_{\a,\b}$ is not transitive on $\PP \setminus\{\a,\b\}$ and has no fixed points in $\PP \setminus\{\a,\b\}$ (since $G_{\a,\b}$ has no orbit of length 1 on $\PP \setminus\{\a,\b\}$ as seen in \cite{FFXZ, Kantor85}). In case (v) with $q=3$, $G_{\{\a,\b\},\s}$ is not transitive on $\PP \setminus\{\a,\b,\s\}$. Therefore, none of these cases produces any $G$-flag graph.

In case (x), $G_{\{\a, \b\}}$ is $2$-transitive on $\PP \setminus\{\a,\b\}$, and thus $10\cdot 9$ divides $2 |M_{11}|/(12\cdot 11)$, which is a contradiction.

In case (xi), $G_{\a,\b}$ has orbit-lengths $1$ (denote this orbit by $\{\s\}$) and $12$ on $\PP \setminus\{\a, \b\}$. If $G_{\{\a, \b\}}\not\leq G_\s$, then $G_{\{\a, \b\}}$ is transitive on $\PP \setminus\{\a, \b\}$ which has 13 elements, a contradiction. Hence $G_{\{\a,\b\}}\leq G_{\s}$, and $G_{\{\a,\b\}, \s}$ is transitive on $\PP \setminus\{\a, \b, \s\}$. It follows that $\Om:=(\s, \PP \setminus\{\a,\b\})^G$ is a feasible $G$-orbit on the flag set of $\overline{\DD}$. We may identify $G$ with a subgroup of $\GL(4,2)\cong A_8$ of index $8$ isomorphic to $A_7$ acting $2$-transitively on $\mathbb{F}_2^4\setminus\{(0,0,0,0)\}$ in its natural action. By \magma \cite{magma} there are four self-paired $G$-orbits on $\C(\overline{\DD}, \Om)$. Moreover, the order of the four corresponding $G$-flag graphs is $|\Om|=105$, the valencies of them are $12$, $12$, $24$ and $24$, respectively, and in each of them there are six vertices in $\Om(\s)$ which have neighbours in $\Om(\t)$ ($\t\neq\s$).

In case (xii), suppose that there exists a feasible $G$-orbit on the flag set of $\overline{\DD}$. When $u=22$, if $G_{\{\a,\b\}}$ is 2-transitive on $\PP \setminus\{\a, \b\}$, then $20\cdot 19$ divides $2 |M_{22}|/(22\cdot 21)$, which is a contradiction. When $u=11, 12, 23, 24$, we have $G=M_u$, $G_{\{\a,\b\}}$ is 2-transitive on $\PP \setminus\{\a, \b\}$ as $G$ is 4-transitive on $\PP$, and thus the set of flags of $\overline{\DD}$ is a feasible $G$-orbit. Denote the set of flags of $\overline{\DD}$ by $\Om$. Since we will take the advantage of permutation representations of $M_{u}$ ($u=11, 12, 23, 24$) from \cite[P.209, Exercise 6.8.10]{Dixon-Mortimer}, we may assume that $\PP=\{1,2,\ldots,u\}$. Since we will use \magma to do computing, for simplicity we replace $0,1,\ldots,u-1$ in \cite[P.209, Exercise 6.8.10]{Dixon-Mortimer} by $1,2,\ldots,u$, respectively. Let $\Psi$ be an $M_u$-orbit on $\C(\overline{\DD}, \Om)$. We may assume $\Psi=((3, \PP\setminus\{1, 2\}), (4, \PP\setminus N))^{M_{u}}$, where $N$ is a $2$-subset of $\{1,2,5,6,\ldots,u-1,u\}$. In the case when $u=24$, $M_{24}=\la a,b,c\ra\leq S_{24}$, where
\begin{align*}
&a=(1,2,3,4,5,6,7,8,9,10,11,12,13,14,15,16,17,18,19,20,21,22,23),\\ &b=(3,17,10,7,9)(5,4,13,14,19)(11,12,23,8,18)(21,16,15,20,22),\\ &c=(1,24)(2,23)(3,12)(4,16)(5,18)(6,10)(7,20)(8,14)(9,21)(11,17)(13,22)(19,15).
\end{align*}
Using \magma we obtain that there are four $G$-orbits on $\C(\overline{\DD}, \Om)$ and all of them are self-paired. Moreover, the four corresponding $G$-flag graphs all have order $6072$, there are $231$ vertices in $\Om(3)$ which have neighbours in $\Om(4)$, and the vertex $(3, \PP\setminus\{1,2\})$ is adjacent to $1$, $30$, $40$ and $160$ vertices in $\Om(4)$, respectively. In the case when $u=23$, there are five self-paired (and some non-self-paired) $G$-orbits on $\C(\overline{\DD}, \Om)$. The five corresponding $G$-flag graphs have order $5313$, there are $210$ vertices in $\Om(3)$ which have neighbours in $\Om(4)$, and the vertex $(3, \PP\setminus\{1,2\})$ is adjacent to $1$, $3$, $6$, $24$ and $32$ vertices in $\Om(4)$, respectively. In the case when $u=12$, there are five $G$-orbits on $\C(\overline{\DD}, \Om)$ and each of them is self-paired. In each of the five corresponding $G$-flag graphs the number of vertices is $660$ and there are $45$ vertices in $\Om(3)$ which have neighbours in $\Om(4)$. Moreover, in these five graphs the vertex $(3, \PP\setminus\{1,2\})$ is adjacent to $1$, $4$, $8$, $16$ and $16$ vertices in $\Om(4)$, respectively. Finally, in the case when $u=11$, there are six self-paired (and some non-self-paired) $G$-orbits on $\C(\overline{\DD}, \Om)$. The corresponding $G$-flag graphs have $495$ vertices, with $36$ vertices in $\Om(3)$ having neighbours in $\Om(4)$, and the vertex $(3, \PP\setminus\{1,2\})$ is adjacent to $1$, $1$, $2$, $2$, $2$ and $2$ vertices in $\Om(4)$, respectively.

\subsection{Affine case}
\label{sec:affine simple 2-transitive groups}

Suppose that $G$ is $2$-transitive on $\PP$ with an abelian socle. We may assume $\PP$ is a vector space over a finite field. Then $G\leq\AGammaL(n, q)$ and $\PP=\mathbb{F}_q^n$ for some positive integers $n$ and $q$ with $q^n=p^d$, where $p$ is a prime and $d$ a positive integer. Denote the zero vector of $\PP$ by $\bzero$.

\smallskip
\textsf{Case (A):} $n=1$. Suppose that $\Om = (\s, \PP \setminus \{0, 1\})^G$ is feasible. Then $G_{\{0,1\}}$ has an orbit of length at least $q-3$ on $\PP \setminus\{0, 1\}$. Hence $q-3\leq |G_{\{0,1\}}|\leq 2d$, which implies $(p, d)=(2,2), (2, 3), (3,1)$ or $(5,1)$. By our assumption $|\PP|>3$ and by Lemma \ref{lem:two 2-transitive actions}, the cases $(p, d)=(2, 3)$ and $(p, d)=(3, 1)$ do not produce any $G$-flag graph. In the case when $(p, d)=(2,2)$, we have $\PP=\{0,1,\a,\a^2\}$, where $\a$ is a generator of $\mathbb{F}_4^\times$, and $G_0=\GL(1,4)$ or $\GammaL(1,4)$. The set of flags of $\overline{\DD}$ is the only feasible $G$-orbit, and $\Psi=((\a, \PP \setminus \{0, 1\}),(\a^2, \PP \setminus \{0, 1\}))^G$ is the only self-paired $G$-orbit on $\C(\overline{\DD}, \Om)$. The corresponding $G$-flag graph $\Ga:=\Ga(\overline{\DD}, \Om, \Psi)$ has $12$ vertices and $\Ga[\Om(\a), \Om(\a^2)]\cong K_2$. In the case when $(p, d)=(5,1)$, we have $\PP=\{0,1,2,3,4\}$ and $G_0=\GL(1,5)$. Moreover, $\Om = (3, \PP\setminus\{0, 1\})^G$ is the only feasible $G$-orbit and $\Psi=((3, \PP \setminus \{0, 1\}),(2, \PP \setminus \{0,4\}))^G$ is the only self-paired $G$-orbit on $\C(\overline{\DD}, \Om)$. The corresponding $G$-flag graph $\Ga:=\Ga(\overline{\DD}, \Om, \Psi)$ has $10$ vertices and $\Ga[\Om(0), \Om(1)]\cong K_2$.

\smallskip
\textsf{Case (B):} $n\geq2$ and $q\geq4$. Suppose that $\Om = (\s, \PP \setminus \{\bzero, \bfe_1\})^G$ is feasible. Choose $\g\in \la\bfe_1\ra \setminus \{\bzero, \bfe_1\}$ not fixed by $G_{\{\bzero,\bfe_1\}}$. (Such an element $\g$ exists as $G_{\{\bzero,\bfe_1\}}$ stabilizes at most one point on $\PP\setminus\{\bzero,\bfe_1\}$ by the feasibility of $\Om$.) Then $\PP \setminus \{\bzero, \bfe_1, \s\}\subseteq\g^{G_{\{\bzero,\bfe_1\}, \s}}\subseteq\g^{G_{\{\bzero,\bfe_1\}}}$. On the other hand, we have $\g^{G_{\{\bzero,\bfe_1\}}}=\g^{G_{\bzero,\bfe_1}}\cup (\g^g)^{G_{\bzero,\bfe_1}}$, where $g\in G_{\{\bzero,\bfe_1\}}\setminus G_{\bzero,\bfe_1}$, and it can be verified that $\g^{G_{\{\bzero,\bfe_1\}}}\subseteq\la\bfe_1\ra$. This is a contradiction and thus this case does not produce any $G$-flag graph.


\medskip
\textsf{Case (C):} $n\geq2$, $q=2$ or $3$.

\medskip
\textsf{Subcase (C1):} $q=2$. In the case when $n=2$, suppose that there exists a feasible $G$-orbit $\Om$ on the flag set of $\overline{\DD}$. Then $G=\mathbb{F}_2^2\rtimes\la\bmat{1 & 1\\1 & 0}\ra\cong A_4$ or $G=\AGL(2,2)\cong S_4$, and the set of flags of $\overline{\DD}$ is the only feasible $G$-orbit. Moreover, $\Psi=((\bfe_2, \PP \setminus \{\bzero, \bfe_1\}),(\bfe_1+\bfe_2, \PP \setminus \{\bzero, \bfe_1\}))^G$ is the only self-paired $G$-orbit on $\C(\overline{\DD}, \Om)$. The corresponding $G$-flag graph $\Ga:=\Ga(\overline{\DD}, \Om, \Psi)$ has $12$ vertices and $\Ga[\Om(\bzero), \Om(\bfe_1)]\cong K_2$.

In the rest of this subcase we assume that $n\geq3$. Suppose that $\Om = (\s, \PP \setminus \{\bfe_1, \bfe_2\})^G$ is a feasible $G$-orbit. By Lemma \ref{lem:two 2-transitive actions}, $G_{\{\bfe_1, \bfe_2\}}$ should be $2$-transitive on $\PP\setminus\{\bfe_1, \bfe_2\}$, and thus we may assume $\s=\bzero$. Since $\Om$ is feasible by our assumption, $G_{\bzero, \{\bfe_1, \bfe_2\}}$ is transitive on $\PP\setminus\{\bzero, \bfe_1, \bfe_2\}$. Hence there exists some $g\in G_{\bzero, \{\bfe_1, \bfe_2\}}:\bx\mapsto \bx A$, where $A\in \GL(n, 2)$, such that $\bfe_n^g=\bfe_1+\bfe_2$. This implies $\row_n(A)=\row_1(A)+\row_2(A)$ and thus $A\not\in \GL(n, 2)$, a contradiction. Therefore, we do not get any $G$-flag graph when $q=2$ and $n \ge 3$.

\medskip
\textsf{Subcase (C2):} $q=3$. In this case $G_{\{\bfe_1, -\bfe_1\}}$ stabilizes ${\bzero}$. Thus, by Lemma \ref{lem:at most one feasible orbit}, $\Om = (\bzero, \PP \setminus \{\bfe_1, -\bfe_1\})^G$ is the only possible feasible $G$-orbit on the set of flags of $\overline{\DD}$. Moreover, a flag $(\s, \PP\setminus\{\a, \b\})\in \Om$ if and only if $\a+\b=-\s$. Suppose that $\Om$ is feasible. Let $\Psi$ be a $G$-orbit on $\C(\overline{\DD}, \Om)$. Since $G_{\{\bfe_1, -\bfe_1\}}$ stabilizes $\bzero$ and is transitive on $\PP \setminus\la\bfe_1\ra$, we may assume $\Psi= ((\bzero, \PP\setminus\{\bfe_1, -\bfe_1\}), (\bfe_n, \PP\setminus N))^G$ for some $2$-subset $N$ of $\PP\setminus\{\bzero, \bfe_n\}$. Since $(\bfe_n, \PP\setminus N)\in \Om$, we may assume $N=\{\bx, -\bx-\bfe_n\}$, where $\bx\not \in\{\bzero, \bfe_n, -\bfe_n\}$. If $\Psi$ is self-paired, that is, if there exists $g\in G$ interchanging $(\bzero, \PP\setminus\{\bfe_1, -\bfe_1\})$ and $(\bfe_n, \PP\setminus N)$, then $g=\d(A, \bfe_n):\bx\mapsto \bx A+\bfe_n$ for some $A\in G_{\bzero}$ and the following relations hold:
\begin{equation}
 \label{equ:self-parity condition q=3, (1)}
\bfe_n A+\bfe_n=\bzero,\;\; \{\bfe_1, -\bfe_1\}^g=N,\;\;\{\bfe_1, -\bfe_1\}^{g^2}=\{\bfe_1, -\bfe_1\}.
\end{equation}
Without loss of generality we may assume $\bfe_1 A + \bfe_n = \bx$, and thus $\bfe_1 A^2 = (\bx -\bfe_n)A=\bfe_1\text{\; or \,} -\bfe_1$ as $g^2=\d(A^2, \bzero)$. Hence \eqref{equ:self-parity condition q=3, (1)} is equivalent to the following:
\begin{equation}
\label{equ:self-parity condition q=3, (2)}
\bfe_n A = -\bfe_n,\;\; (\bfe_1 - \bfe_n) A = \bx,\;\;\bx A=\bfe_1 - \bfe_n \text{\; or \,}-\bfe_1 - \bfe_n.
\end{equation}
Therefore, a $G$-orbit $\Psi$ on $\C(\overline{\DD}, \Om)$ is self-paired if and only if there exists some $A\in G_\bzero$ such that \eqref{equ:self-parity condition q=3, (1)} or \eqref{equ:self-parity condition q=3, (2)} holds.

The affine $2$-transitive permutation groups $G$ are all known. The degree $u:=|\PP|= p^d$ and the stabilizer $G_{\bf 0}$ in $G$ of the zero vector ${\bf 0}$ are as follows (\cite{Kantor85}, \cite[p.194]{Cameron99}, \cite{Cameron81}, \cite[p.386]{Huppert82}), where $q = p = 3$ and $n=d$ in all cases:
\begin{itemize}[topsep=0.5ex]
\setlength\itemsep{-0.3em}
\item[(i)] $G_{\bf 0} \leq \GammaL(1, q)$, $q=p^d$;
\item[(ii)] $G_{\bf 0} \unrhd \SL(n, q)$, $n \geq 2$, $q^n = p^d$;
\item[(iii)] $G_{\bf 0} \unrhd \Sp(n, q)$, $n \geq 4$, $n$ is even, $q^n = p^d$;
\item[(iv)]  $G_{\bf 0} \unrhd G_2(q)$, $q^6 = p^d$, $q >2$, $q$ is even;
\item[(v)] $G_{\bf 0} = G_2(2)' \cong \PSU(3, 3)$, $u = 2^6$;
\item[(vi)] $G_{\bf 0} \cong A_6$ or $A_7$, $u = 2^4$;
\item[(vii)] $G_{\bf 0} \unrhd \SL(2, 5)$ or $G_{\bf 0} \unrhd \SL(2, 3)$, $d=2$, $p = 5, 7, 11, 19, 23, 29$ or $59$;
\item[(viii)] $G_{\bf 0} \cong \SL(2,13)$, $u = 3^6$;
\item[(ix)] $d=4$, $p=3$, $G_{\bf 0} \unrhd \SL(2, 5)$ or $G_{\bf 0} \unrhd E$, where $E$ is an extraspecial group of order $32$.
\end{itemize}

\delete
{
\begin{description}
\item[\rm (i)] $G_{\bf 0} \leq \GammaL(1, q)$, $q=p^d$;

\item[\rm (ii)] $G_{\bf 0} \unrhd \SL(n, q)$, $n \geq 2$, $q^n = p^d$;

\item[\rm (iii)] $G_{\bf 0} \unrhd \Sp(n, q)$, $n \geq 4$, $n$ is even, $q^n = p^d$;

\item[\rm (iv)]  $G_{\bf 0} \unrhd G_2(q)$, $q^6 = p^d$, $q >2$, $q$ is even;

\item[\rm (v)] $G_{\bf 0} = G_2(2)' \cong \PSU(3, 3)$, $u = 2^6$;

\item[\rm (vi)] $G_{\bf 0} \cong A_6$ or $A_7$, $u = 2^4$;

\item[\rm (vii)] $G_{\bf 0} \unrhd \SL(2, 5)$ or $G_{\bf 0} \unrhd \SL(2, 3)$, $d=2$, $p = 5, 7, 11, 19, 23, 29$ or $59$;

\item[\rm (viii)] $G_{\bf 0} \cong \SL(2,13)$, $u = 3^6$;

\item[\rm (ix)] $d=4$, $p=3$, $G_{\bf 0} \unrhd \SL(2, 5)$ or $G_{\bf 0} \unrhd E$, where $E$ is an extraspecial group of order $32$.
\end{description}
}




In case (i), the cases $(p, d)=(2,2)$ and $(p, d)=(5,1)$ can occur and have been considered in Case (A) above.

In case (ii), we only need to consider the case when $q=3$. First assume $n=2$. Then $G_0=\SL(2,3)$ or $\GL(2,3)$. Let $\Om$ and $\Psi$ be as in Subcase (C2) above ($\Om$ is feasible). When $G_0=\SL(2,3)$, by \eqref{equ:self-parity condition q=3, (2)}, $\Psi$ is self-paired if and only if $N=\{-\bfe_1+\bfe_2, \bfe_1+\bfe_2\}$. The corresponding $G$-flag graph $\Ga$ has $36$ vertices, and $\Ga[\Om(\bzero), \Om(\bfe_2)]\cong 3\cdot K_2$. When $G_0=\GL(2,3)$, by \eqref{equ:self-parity condition q=3, (2)}, $\Psi$ is self-paired if and only if $N=\{\bfe_1, -\bfe_1-\bfe_2\}$, $\{-\bfe_1, \bfe_1-\bfe_2\}$ or $\{-\bfe_1+\bfe_2, \bfe_1+\bfe_2\}$. Since $((\bzero, \PP\setminus\{\bfe_1, -\bfe_1\}), (\bfe_2, \PP\setminus \{\bfe_1, -\bfe_1-\bfe_2\}))$ and $((\bzero, \PP\setminus\{\bfe_1, -\bfe_1\}), (\bfe_2, \PP\setminus \{-\bfe_1, \bfe_1-\bfe_2\}))$ are in the same $G$-orbit, there are exactly two self-paired $G$-orbits on $\C(\overline{\DD}, \Om)$. Moreover, the two corresponding $G$-flag graphs $\Ga$ both have $36$ vertices, the number of vertices in $\Om(\bzero)$ which have neighbours in $\Om(\bfe_2)$ is $3$, and the valencies of $\Ga[\Om(\bzero), \Om(\bfe_2)]$ are $2$ and $1$, respectively. Next assume $n\geq3$. Then $G_0=\SL(n,3)$ or $\GL(n,3)$. Let $\Om$ and $\Psi$ be as in Subcase (C2) above (we can see that $\Om$ is feasible). When $\bx\in \la \bfe_1,\bfe_n\ra$, $\Psi$ is self-paired if and only if $N=\{\bfe_1, -\bfe_1-\bfe_n\}$, $\{-\bfe_1, \bfe_1-\bfe_n\}$ or $\{-\bfe_1+\bfe_n, \bfe_1+\bfe_n\}$. Since $((\bzero, \PP\setminus\{\bfe_1, -\bfe_1\}), (\bfe_n, \PP\setminus \{\bfe_1, -\bfe_1-\bfe_n\}))$ and $((\bzero, \PP\setminus\{\bfe_1, -\bfe_1\}), (\bfe_n, \PP\setminus \{-\bfe_1, \bfe_1-\bfe_n\}))$ are in the same $G$-orbit, they produce the same $G$-flag graph. Moreover, the corresponding $G$-flag graphs $\Ga$ have $3^n(3^{n}-1)/2$ vertices, the number of vertices in $\Om(\bzero)$ which have neighbours in $\Om(\bfe_n)$ is $(3^n-3)/2$, and the valencies of $\Ga[\Om(\bzero), \Om(\bfe_n)]$ are $2$ and $1$, respectively. When $\bx\not\in \la \bfe_1,\bfe_n\ra$, it can be verified that $\Psi$ is self-paired, and all $\bx\not\in \la \bfe_1,\bfe_n\ra$ produce the same $\Psi$ and thus the same $G$-flag graph. The corresponding $G$-flag graph $\Ga$ have $3^n(3^{n}-1)/2$ vertices, the number of vertices in $\Om(\bzero)$ which have neighbours in $\Om(\bfe_n)$ is $(3^n-3)/2$, and the valency of $\Ga[\Om(\bzero), \Om(\bfe_n)]$ is $(3^n-9)/2$.

In case (iii), we only need to consider the case when $q=3$. We may assume that $\PP$ is equipped with a nonsingular alternating bilinear form $B(\bx, \by)$:
$$
B(\bx, \by)=\sum_{i=1}^{n} \left(x_{2i-1} y_{2i} - x_{2i} y_{2i-1}\right),
$$
where $\bx=(x_1,x_2,\ldots,x_{2n-1}, x_{2n})$, $\by=(y_1, y_2,\ldots,y_{2n-1}, y_{2n})\in\PP$. As in \cite{Dye}, denote by $\GSp(n,3)$ the subgroup consisting of all $A\in\GL(n,3)$ satisfying $B(\bx A, \by A)=\l B(\bx,\by)$ for any $\bx,\by\in \PP$, where $\l=1$ or $-1$. Then $[\GSp(n,3):\Sp(n,3)]=2$. Let $L$ be the normalizer of $\Sp(n,3)$ in $\GL(n,3)$. If $L>\GSp(n,3)$, then $L\geq\SL(n,3)$ by the proof of \cite[Theorem 3]{Dye} and thus $\Sp(n,3)\unlhd\SL(n,3)$, which is a contradiction as $\PSL(n,3)$ is simple. Hence the normalizer of $\Sp(n,3)$ in $\GL(n,3)$ is $\GSp(n,3)$ and we obtain $G_\bzero=\Sp(n,3)$ or $\GSp(n,3)$. Let $\Om$ be as in Subcase (C2) above. Suppose that $\Om$ is feasible. Then $G_{\bzero,\{\pm\bfe_1\}}$ is transitive on $\PP\setminus\{\bzero,\pm\bfe_1\}$ and there exists $A\in G_\bzero$ such that $\bfe_1 A=\pm\bfe_1$ and $\bfe_n A=\bfe_2$. This implies $\pm1=B(\pm\bfe_1,\bfe_2)=B(\bfe_1 A, \bfe_n A)=\pm B(\bfe_1, \bfe_n)=0$, a contradiction. Therefore, case (iii) does not produce any $G$-flag graph.

Cases (iv)--(vii) do not produce any $G$-flag graph based on the above discussion.

In case (viii), since $|G_{\{\bzero,\bfe_1\}}|=6$, $G_{\{\bzero,\bfe_1\}}$ cannot be transitive on $\PP \setminus\{\bzero, \bfe_1, -\bfe_1\}$ which has $3^6-3=726$ elements. This case does not produce any $G$-flag graph.

In case (ix), if $G_{\bf 0} \unrhd \SL(2, 5)$, then $|G_{\{\bzero,\bfe_1\}}|\leq2\cdot12$ (see \cite{FFXZ2}), and $G_{\{\bzero,\bfe_1\}}$ cannot be transitive on $\PP \setminus\{\bzero, \bfe_1, -\bfe_1\}$ which has $3^4-3=78$ elements. If $G_{\bf 0} \unrhd E$, where $E$ is an extraspecial group of order $32$, then $|G_{\{\bzero,\bfe_1\}}|=2 |G_{\bzero, \bfe_1}|$ divides $96$ (see \cite{Kantor85}, \cite{FFXZ2}), and $G_{\{\bzero,\bfe_1\}}$ cannot be transitive on $\PP \setminus\{\bzero, \bfe_1, -\bfe_1\}$ which has $3^4-3=78$ elements. Hence this case does not produce any $G$-flag graph.

So far we have completed the proof of part (b) of Theorem \ref{thm:main theorem} and hence the proof of the whole Theorem \ref{thm:main theorem}.




\medskip
\medskip
\noindent\textbf{Acknowledgement}~~Shenglin Zhou was supported by the National Natural Science Foundation of China (Grant No.12271173).

\medskip
\medskip
\noindent\textbf{Statements and Declarations}~~The authors have no relevant financial or non-financial interests to disclose.

\newpage

\appendix
\section*{Appendix: Sample \magma codes (for referees only)}

The following \magma codes are for case (xii) in Section \ref{sec:almost simple 2-transitive groups}. For other cases where we use \magma to do the calculations, the \magma codes are similar.

\begin{verbatim}
S:=Sym(24);
a:=S!(1,2,3,4,5,6,7,8,9,10,11,12,13,14,15,16,17,18,19,20,21,22,23);
b:=S!(3,17,10,7,9)(5,4,13,14,19)(11,12,23,8,18)(21,16,15,20,22);
c:=S!(1,24)(2,23)(3,12)(4,16)(5,18)(6,10)(7,20)(8,14)(9,21)(11,17)(13,22)
(19,15);
G:=sub<S|a,b,c>;               C:={1,2};                  D:={3,4};
G3:=Stabilizer(G,3);           J:=Stabilizer(G3,C);       #G3/#J;
G34:=Stabilizer(G3,4);         J2:=Stabilizer(G34,C);     #G34/#J2;
GD:=Stabilizer(G,D);           X:=[g:g in GD|3^g eq 4];
Y:=[g:g in X|(C^g)^g eq C];    P:=[C^(Y[1])];

for g in Y do
  B:=C^g;i:=#P;
    for j in [1..i] do
        if B eq P[j] then break j;
          else if j eq i then P:=Append(P,B);
               end if;
        end if;
    end for;
end for;

#P;

GC:=Stabilizer(G,C);    Z:=Stabilizer(GC,3);   H:=Stabilizer(Z,4);
h:=#H;       Q:=[];      x:=[];       y:=[];

for N in P do
  m:=#Stabilizer(H,N);n:=h/m;i:=#Q;
   if i eq 0 then Q:=Append(Q,n); x:=Append(x,1); y:=Append(y,N); else
    for j in [1..i] do
     if n eq Q[j] then x[j]:=x[j]+1; break j;
      else if j eq i then Q:=Append(Q,n); x:=Append(x,1); y:=Append(y,N);
       end if;
     end if;
    end for;
   end if;
end for;

#Q;#x;#y;    Q;x;y;
\end{verbatim}

\end{document}